\magnification=\magstep1
\everymath={\displaystyle}

\parskip=\medskipamount
%
\centerline {\bf The $q$-Harmonic Oscillator and an Analog
                  of the Charlier polynomials}
\bigskip
\centerline {R. Askey \plainfootnote{\dag}{Department of
Mathematics,
University of Wisconsin, Madison, WI 53706, USA}
and S. K. Suslov \plainfootnote{\ddag}{Russian Scientific Center
\lq \lq
Kurchatov Institute", Moscow 123182, Russia}}
\medskip
{\parindent = 12pt \narrower\narrower\narrower
\noindent{\bf Abstract}. A model of a $q$-harmonic oscillator based
on
$q$-Charlier polynomials of Al-Salam and Carlitz is discussed.
Simple explicit realization of $q$-creation and $q$-annihilation
operators, $q$-coherent states and an analog of the Fourier
transformation
are found. A connection of the kernel of this transform with
biorthogonal
rational functions is observed. \par}
\medskip
Models of $q$-harmonic oscillators are being developed in
connection with
quantum groups and their various applications ( see, for example,
Refs.
[1--5]). The $q$-analogs of boson operators were introduced
explicitly
in Refs. [1,3] and [5], where the corresponding wave functions were
found
in terms of the continuous $q$-Hermite polynomials of Rogers [6,7]
and
in terms of  the Stieltjes--Wigert polynomials [8,9], respectively.
Here we introduce one more explicit realization of $q$-creation and
$q$-annihilation operators with the aid of $q$-Charlier polynomials
of
Al-Salam and Carlitz [10].
\medskip
The $q$-orthogonal polynomials $V_n^{a}$ studied by Al-Salam and
Carlitz
may be considered as a $q$-{\it version of the Charlier
polynomials\/}
$c_n^{\mu}(s)$
( see, for example, [11,12] ). To emphasize this analogy we use the
notation
$c_n^{\mu}(x\mid q)$ for the Al-Salam and Carlitz polynomials. In
our
notation  they can be defined by the three-term recurrence relation
$$
\mu q^{-n-1}c_{n+1}^\mu (x\mid q) +(1-q^n)q^{-n}c_{n-1}^{\mu
}(x\mid q)=
\left( (\mu +q)q^{-n-1}-x\right) c_n^{\mu}(x\mid q)\,, \eqno (1)
$$
$c_0^{\mu}(x\mid q)=1\,,c_1^{\mu}(x\mid q)=\mu^{-1}(\mu +q-qx)\,.$
These
polynomials are orthogonal
$$\sum_{s=0}^{\infty}c_m^{\mu}(q^{-s}\mid q)\,c_n^{\mu}(q^{-s}\mid
q)\,
\rho(s)q^{-s}={(q;q)_n \over \mu^n }\delta _{mn} \eqno (2)
$$
with respect to a positive measure
$$
\rho(s)=(\mu ;q)_{\infty}{\mu^sq^{s^{2}} \over (q,\mu ;q)_s};\,
0<\mu ,q<1\,,
\eqno (3)
$$
where the usual notations (see [13]) are
$$
\eqalignno{
&(a;q)_n=\prod_{k=0}^{n-1}(1-aq^k)\,,\cr
&(a,b;q)_n=(a;q)_n(b;q)_n\,,&(4)\cr
&(a;q)_{\infty}=\lim_{n\to \infty }(a;q)_n\,.\cr }
$$
The weight function (3) is a solution  of the Pearson equation
$\Delta
(\sigma \rho) = \rho \tau \nabla x_1 \;$(for details,
see [12,14]) with\;$ x (s) = q^{-s},\; \sigma (s) =(1-q^{-s})(\mu
-q^{1-s}),\;
\sigma (s) + \tau (s) \nabla x_1 (s) = \mu$. The explicit form of
the polynomials $c_n^{\mu }(x\mid q)$ is
$$
c_n^{\mu }(x\mid q) = \,_2\varphi_0(q^{-n},x;q,\,q^{1+n}/\mu
),\,x=q^{-s}\,.
\eqno (5)
$$
For the definition of the basic hypergeometric function $_2 \varphi
_0 $\,,
see [13].
In the limit $q\to 1$ it easy to see from (1) or (5) that
$$
\lim _{q\to 1}c_n^{(1-q)\mu} (q^{-s}\mid q) = \,_2F_0 (-n,-s;-1/\mu
) =
c_n^{\mu}(s)\,. \eqno (6)
$$
This justifies our notation for the Al-Salam and Carlitz 
polynomials.
\medskip
The polynomials $c_n^{\mu }(x\mid q)$ give us the possibility to
introduce
a new model of a $q$-oscillator. We can define a $q$-{\it version
of the wave
functions\/} of harmonic oscillator as
$$
\psi_n(s)=d^{-1}_n q^{-{s\over2}} \rho^{1 \over 2} (s)\, c_n^{\mu
}(q^{-s}
\mid q) \,, \eqno (7)
$$
where $d^2_n=(q;q)_n/\mu ^n$. These $q$-wave functions satisfy the
orthogonality relation
$$
\sum^\infty _{s=0} \psi_n(s)\,\psi_m(s)= \delta_{nm} \,. \eqno (8)
$$
The $q$-{\it annihilation\/} $a$ and $q$-{\it creation\/} $a^+$
operators
have the following explicit form
$$\eqalignno{
a=&(1-q)^{- {1 \over 2}} \left[\mu ^{1 \over 2} q^s -
\sqrt{(1-q^{s+1})(1-\mu q^s)}
e^{\partial _s} \right]\,,\cr
&&(9) \cr
a^+=&(1-q)^{- {1 \over 2}} \left[\mu ^{1 \over 2} q^s -
e^{-\partial _s}
 \sqrt {(1-q^{s+1}) (1-\mu q^s)} \right]\,,\cr
}$$
where $\partial _s \equiv {d \over ds},\; e^{\alpha \partial_s}
f(s)=f(s+
\alpha)$.  They satisfy the $q$-commutational rule
$$
aa^+-qa^+a=1\eqno {(10)}
$$
and act on the $q$-wave functions defined in (7) by
$$
a\psi_n=e^{1 \over 2}_{n\vphantom1}\psi_{n-1},\;
a^+\psi_n=e^{1 \over2}_{n+1}\psi_{n+1}\,,\eqno (11)
$$
where
$$
e_n={1-q^n \over 1-q}\,.
$$
In this model of the $q$-oscillator equations (11) are equivalent
to
difference-differentiation formulas
$$\eqalign{
&\mu q^s \Delta c_n^{\mu}(x \mid q)=(q^n-1)\,c_{n-1}^{\mu}(x \mid
q)\,,\cr
&q^s \nabla [ \rho(s) c_n^{\mu}(x \mid q)]= \rho(s)
\,c_{n+1}^{\mu}(x \mid q)\,,
\cr}$$
respectively. Here $\Delta f(s)= \nabla f(s+1)=f(s+1)-f(s)$ and
$x=q^{-s}$.
In view of (6) the functions $\psi _n(s)$ converge in the limit
$q\to 1^-$
to the wave functions of the discrete model of the linear harmonic
oscillator
considered in [15].
\medskip
The $q$-{\it Hamiltonian\/} $H=a^+a$ acts on the wave functions (7)
as
$$
H\psi_n=e_n\psi_n\,\eqno(12)
$$
and has the following explicit form
$$\eqalignno{
H &= (1-q)^{-1} \bigl[\mu q^{2s}+(1-q^s)(1- \mu q^{s-1} )- &(13)
\cr
 & \mu ^ {1 \over 2}q^s \sqrt{(1-q^{s+1})(1-\mu q^s)}e^{\partial
_s}-
\mu ^ {1 \over 2}q^{s-1} \sqrt{(1-q^s)(1-\mu q^{s-1}}e^{-\partial
_s} \bigr]
\,. \cr }
$$
By factorizing the Hamiltonian ( or the difference equation for the
Al-Salam
and Carlitz polynomials ) we arrive at the explicit form (9) for
the $q$-boson
operators.

Since $a^+a=H$,  the relation (10) can be written in the equivalent
form
$$
[a,a^+]=1-(1-q)H\equiv q^N\,. \eqno(14)
$$
The operator
$$
N={1 \over {\log q}}\log [1-(1-q)H]\eqno {(15)}
$$
can be considered as the {\it number operator\/}, since
$$
[a,N]=a,\;\qquad  [N,a^+]=a^+\,.\eqno{(16)}
$$
From these relations one can obtain the equations (11) and the
spectrum (12)
of the $q$-Hamiltonian in abstract form. The $q$-wave functions are
$$
\psi _n(s)=c_n \left(a^+ \right)^n \psi _0(s)\,,\qquad a\, \psi
_0(s)=0\,,
$$
where $c_n=(e_n!)^{-1/2}$ and $ e_n!=e_1e_2 \dots e_n $.
\medskip
For the model of the $q$-oscillator under discussion we can
construct
explicitly $q$-{\it coherent states\/} and an analog of the Fourier
transformation.
For the coherent states $\mid \alpha \rangle$ defined by
$$
\eqalignno{
&a\mid \alpha \rangle =\alpha \mid \alpha \rangle \,, \qquad
\langle \alpha |\alpha \rangle =1\,,\cr
&\mid \alpha \rangle
=f_{\alpha}\sum_{n=0}^{\infty}{\alpha^n\psi_n(s)\over
(e_n!)^{1/2}}\,,&(17)\cr
&f_{\alpha}=\left((1-q)\mid \alpha |^2;q\right)^{1\over 2}_{ \infty
}\,,
\quad (1-q)\mid \alpha |^2 < 1 \cr}
$$
we can write
$$
\mid \alpha \rangle = f_{\alpha} \left( \rho q^{-s} \right)^{1
\over 2}
\sum_{n=0}^\infty
{t^n \over (q;q)_n}\,c_n^{\mu}(x \mid q)\,,t=\alpha \mu^{1 \over 2}
(1-q)^{1 \over 2}\,. \eqno{(18)}
$$
With the aid of the generating function [10]
$$
\sum_{n=0}^\infty {t^n \over (q;q)_n}\,c_n^{\mu}(x \mid q)=
{(qtx/ \mu ;q)_\infty \over (t\,,qt/ \mu ;q)_\infty }\,,\,
|t|<1,\, |qt/\mu |<1
\eqno{(19)}
$$
we arrive at the following {\it explicit form \/} for the
$q$-coherent states
$$
\mid \alpha \rangle =f_{\alpha}\left( \rho q^{-s} \right)^{1 \over
2}
{\left(\alpha (1-q)^{1/2} \mu ^{-1/2}q^{1-s};q \right)_\infty
\over \left(\alpha (1-q)^{1/2}\mu ^{1/2},\alpha
q(1-q)^{1/2}\mu^{-1/2};q
\right)_\infty} \,, \eqno{(20)}
$$
where $\rho =(q;q)_\infty^{-1}(q^{s+1}\,,\mu q^s;q)_\infty \mu^s
q^{s^2}\,.$
These coherent states are not orthogonal
$$
\langle \alpha \mid \beta \rangle ={\left((1-q) \mid \alpha \mid
^2\,,
(1-q) \mid \beta
\mid ^2;q \right)_\infty ^{1/2} \over \left((1-q)\alpha^*\beta\,;q
\right)_\infty} \,,
$$
where $*$ denotes the complex conjugate.
\medskip
To define an analog of the {\it Fourier transform \/} we can
consider,
following Wiener's approach to the classical Fourier transform [16]
( see also [17,18] ), the kernel of the form
$$\eqalignno{
K_t(s,p)&=\sum^\infty _{n=0}t^n \psi_n(s)\,\psi_n(p) \cr
&&(21)\cr
&=\left( \rho(s) \rho(p) q^{-s-p} \right)^{1 \over 2}
\sum_{n=0}^\infty {(\mu t)^n \over (q;q)_n}\,
c_n^{\mu}(q^{-s}\mid q)\,c_n^{\mu}(q^{-p}\mid q)\,.\cr }
$$
The series can be summed with the aid of the bilinear generating
function
by Al-Salam and Carlitz [10]
$$
\sum_{n=0}^\infty c_n^{\mu _1}(x\mid q)\,c_n^{\mu _2}(y \mid q)\,
{t^n \over (q;q)_n}=
{(qtx/ \mu_1\,,qty/\mu_2 ;q)_\infty \over
(t\,,qt/ \mu_1\,,qt/\mu_2 ;q)_\infty }
$$
$$
\cdot \,_3\varphi_2\left(
\matrix
x\,, y\,, t\, \\
\noalign{\smallskip}
qtx/\mu_1,  qty/\mu_2 \\
\endmatrix
; q,{q^2 t \over \mu_1 \mu_2} \right) \,.\eqno(22)
$$
The answer is
$$
K_t(s,p)=\left( \rho(s) \rho(p) q^{-s-p} \right)^{1 \over 2}
{(tq^{1-s}\,,tq^{1-p}\,;q)_\infty \over
(qt\,,qt\,,\mu t ;q)_\infty }
$$
$$
\cdot \, _3\varphi_2\left(
\matrix
q^{-s}\,, q^{-p}\,,\mu t\, \\
\noalign{\smallskip}
tq^{1-s}\,, tq^{1-p} \\
\endmatrix
; q,{q^2 t \over \mu} \right) \,. \eqno(23)
$$
The $q$-wave functions (7) are eigenfunctions of the \lq \lq
discrete
$q$-Fourier transform'',
$$
i^m\psi_m(s)=\sum_{p=0}^\infty K_i(s,p)\,\psi_m(p)\,.\eqno(24)
$$
The orthogonality relation of the kernel,
$$
\sum^\infty _{p=0} K_i(s,p)\,K_i^*(s',p)= \delta_{ss'} \,, \eqno
(25)
$$
implies the orthogonality of the rational functions (23). In view
of (6)
in the limit $q \to 1^-$ we get the \lq \lq discrete Fourier
transform''
considered in [17].
\medskip
Similarly, with the aid of the bilinear generating function (22)
and
the orthogonality property of the Wall polynomials, which are dual
to
the polynomials (5), one can obtain the {\it biorthogonality
relation\/}
$$
\sum^\infty _{s=0} u_m(s)\,v_n(s)\, \rho(s)q^{-s}=d^2_n 
\delta_{mn} \eqno (26)
$$
with
$$
\rho(s)={\left({\mu_2 \over t_1 }, {\mu_2 \over t_2 }\,;q \right)_s
\over
\left(q, \mu_2 \,;q \right)_s}\, \mu_1^s q^s
$$
and
$$
d^2_n={\left( t_1,t_2\,;q \right)_\infty \over
\left( \mu_1 ,\mu_2 \,;q \right)_\infty } \cdot
{\left(q, \mu_1 \,;q \right)_n \over \left( {\mu_1 \over t_1 },
{\mu_1 \over t_2 }\,;q \right)_n } \, \mu_2^{-n}
$$
for the $_3 \varphi _2 $-{\it rational functions\/} of the form
$$
u_m(s)= \,_3\varphi_2\left(
\matrix
q^{-m}\,,\, q^{-s}\,, \,t_1 \, \\
\noalign{\smallskip}
{ t_1 \over \mu_1 }q^{1-m}\,,{ t_1 \over \mu_2 } q^{1-s} \\
\endmatrix
; q,\,{q^2 \over t_2 } \right) \,;\eqno(27)
$$
$$
v_n(s)=u_n(s) \vert _{ t_1 \leftrightarrow t_2 }\,; \qquad
t_1 t_2= \mu_1 \mu_2\,.
$$
These functions are self-dual. They belong to  classical
biorthogonal
rational functions [19,20].
\medskip
We have considered here the explicit form of $q$-boson operators
which
satisfy the commutational rule (10) when $0<q<1$. The case $q>1$ is
also
interesting. It leads to another family of Al-Salam and Carlitz
polynomials.
\medskip
It is evident that for the models of the $q$-oscillator under
discussion
one can readily construct dynamical symmetry group $SU_q(1,1)$ [4]
and
write an explicit realization for irreducible representations
$|j,m \rangle _q=\psi_{j+m}(s)\psi_{j-m}(s')$ of the group$SU_q(2)$
[1,2].
\vfill
\eject
\frenchspacing
REFERENCES

\item{1.}  Macfarlane, A.J. (1989) J. Phys. A: Math. Gen., Vol. 22,
p. 4581--4588.

\item{2.}  Biedenharn, L.C. (1989) J. Phys. A: Math. Gen., Vol. 22,
p. L873--L878.

\item{3.}  Atakishiyev, N.M. and Suslov, S.K. (1990) Teor. i Matem.
Fiz.,
Vol. 85, No. 1, p. 64--73.

\item{4.}  Kulish, P.P. and Damaskinsky, E.V. (1990) J. Phys. A:
Math. Gen.,
 Vol. 23, p. L415--L419.

\item{5.}  Atakishiyev, N.M. and Suslov, S.K. (1991) Teor. i Matem.
Fiz.,
Vol. 87, No. 1, p. 154--156.

\item{6.} Rogers, L.J. (1894) Proc. London Math. Soc., Vol. 25, p.
318--343.

\item{7.} Askey, R. and Ismail, M.E.H. (1983) In:{\it  Studies in
Pure
Mathematics\/}
( P. Erd\"os, ed.), Birkh\"auser, Boston, Massachusetts, p. 55--78.

\item{8.}  Stieltjes, T.J. {\it Recherches sur les Fractions
Continues\/},
Annales de la Facult\'e des Sciences de Toulouse,
\underbar {8} (1894) 122 pp., \underbar {9} (1895), 47 pp.
Reprinted in {\it Oeuvres Compl\'etes\/}, vol. 2.

\item{9.}  Wigert S. (1923) Arkiv f\"or Matematik, Astronomi och
Fysik,
Bd. 17, No. 18, p. 1--15.

\item{10.} Al-Salam, W.A. and Carlitz, L. (1965) Math. Nachr., Bd.
30,
S. 47--61.

\item{11.} Chihara, T.S. (1978) {\it An Introduction to Orthogonal
Polynomials.\/} Gordon and Breach, New York.

\item{12.} Nikiforov, A.F., Suslov, S.K., and Uvarov, V.B. (1991)
{\it Classical Orthogonal Polynomials of a Discrete Variable.\/}
Springer-Verlag, Berlin, Heidelberg.

\item{13.} Gasper, G. and Rahman, M. (1990) {\it Basic
Hypergeometric Series,
\/} Cambridge Univ. Press, Cambridge.

\item{14.}  Suslov, S.K. (1989) Russian Math. Surveys, London Math.
Soc.,
Vol. 44, No. 2, p. 227--278.

\item{15.}  Atakishiyev, N.M. and Suslov, S.K. (1989) A Model of
the Harmonic Oscillator on the Lattice. In: {\it Contemporary Group
Analysis:
Methods and Applications,\/} Baku, p. 17--21 [ in Russian ].

\item{16.} Wiener, N. (1933) {\it The Fourier Integral and Certain
of Its
Applications.\/} Cambridge University Press, Cambridge.

\item{17.} Askey, R., Atakishiyev, N.M., and Suslov, S.K.
Fourier Transformations for Difference Analogs of the Harmonic
Oscillator,
to appear

\item{18.} Askey, R., Atakishiyev, N.M., and Suslov, S.K.
An Analog of the Fourier Transformations for a $q$-Harmonic
Oscillator.
Preprint No.5611/1, Kurchatov Institute, Moscow, 1993.

\item{19.} Wilson, J.A. (1991) SIAM J. Math. Anal., Vol. 22(4), p.
1147--1155.

\item{20.} Rahman, M. and Suslov, S.K. Classical Biorthogonal
Rational
Functions, submitted.

\bye